\documentclass{amsart}
\usepackage{amssymb,latexsym,amsmath,amsthm}
\usepackage[mathscr]{eucal}

\newtheorem{theorem}{Theorem}[section]

\numberwithin{equation}{section}

\renewcommand{\l}{\lambda}

\newcommand{\cd}{\,\cdot\,}
\newcommand{\g}{\mathrm{g}}
\newcommand{\dsq}{d^2\!}
\newcommand{\R}{\ensuremath{\mathbb{R}}}

\newcommand{\prtl}{\ensuremath{\partial}}

\newcommand{\hf}{\ensuremath{\frac{1}{2}}}

\title[Multlinear spectral cluster estimates]{On multilinear spectral cluster estimates for manifolds with boundary}

\thanks{The authors were supported by the National Science Foundation
grants DMS-0140499, DMS-0354668, DMS-0555162, and DMS-0354386.}

\author{Matthew D. Blair}
\address{Department of Mathematics, Johns Hopkins University,
Baltimore, MD 21218}
\email{mblair@math.jhu.edu}

\author{Hart F. Smith}
\address{Department of Mathematics, University of Washington,
Seattle, WA 98195}
\email{hart@math.washington.edu}

\author{Christopher D. Sogge}
\address{Department of Mathematics, Johns Hopkins University,
Baltimore, MD 21218}
\email{sogge@jhu.edu}

\begin{document}

\maketitle

\section{Introduction}
Let $(M^n,\g)$ be a smooth, compact $n$-dimensional Riemannian
manifold with boundary and let $\Delta$ be the corresponding
Laplace-Beltrami operator acting on functions.
If the boundary is non-empty, we assume that either
Dirichlet or Neumann conditions are imposed along $\partial M^n$.

Consider the operators $\chi_\lambda$ defined as
projection onto the subspace spanned by the Dirichlet (or Neumann)
eigenfunctions whose corresponding eigenvalues $-\lambda_j^2$
satisfy $\lambda_j \in [\lambda-1, \lambda]$.  In the case
that $\prtl M^n$ is empty,
it was established in~\cite{Sogge88} that the
following, best possible $L^2 \to L^q$ 
estimates hold for $\chi_\lambda$:
\begin{equation}\label{E:nobdcluster}
\|\chi_\lambda\|_{L^2 \to L^q} \lesssim \begin{cases}
\lambda^{\frac{n-1}{2}(\hf-\frac{1}{q})} &  2 \leq q \leq \frac{2(n+1)}{n-1}\\
\lambda^{n(\hf-\frac{1}{q})-\hf} &  \frac{2(n+1)}{n-1} \leq q \leq
\infty
\end{cases}
\end{equation}

Recently, in~\cite{bgtbilin} and~\cite{bgtmultilin}, Burq,
G\'{e}rard, and Tzvetkov established multilinear
versions of these estimates, also under the
assumption that the boundary of $M$ is empty.
To state these, suppose that $\l\ge\mu\ge\nu\ge 1.$ Then
\begin{equation}\label{bilin}
\| \chi_\l f\, \chi_\mu g\|_{L^2(M)} \lesssim 
\Lambda(\mu) \|f\|_{L^2(M)} \|g\|_{L^2(M)}
\end{equation}
\begin{equation}\label{trilin}
\| \chi_\l f \,\chi_\mu g \,\chi_\nu h\|_{L^2(M)} \lesssim
(\mu\nu)^{\frac{2n-3}{4}}\|f\|_{L^2(M)}
\|g\|_{L^2(M)} \|h\|_{L^2(M)}
\end{equation}
where in the first estimate
\begin{equation*}
\Lambda(\mu) = \begin{cases} \mu^\frac{1}{4} & n=2\\
\mu^{\frac{1}{2}}(\log{\mu})^\hf & n=3\\
\mu^{\frac{n-2}{2}} & n\geq 4
\end{cases}
\end{equation*}
With the exception of the logarithmic loss for $n=3$, 
the linear estimate \eqref{E:nobdcluster} with $q=4$ 
follows as a corollary of the bilinear
estimate \eqref{bilin},
by taking $\l=\mu$ and $f=g$.
Similarly, the trilinear estimate \eqref{trilin}
implies \eqref{E:nobdcluster} with $q=6$.
Moreover, by taking $h$ constant and $\nu=1$,
\eqref{trilin} implies \eqref{bilin} in case $n=2$.
For $n\ge 4$, however, the trilinear estimate
can be improved by using \eqref{bilin} together with
the $L^\infty$ bounds \eqref{E:nobdcluster} on $h$.

In the case where $\prtl M^n$ is nonempty, the issue of spectral
cluster estimates is considerably more intricate.  Here the
Rayleigh whispering gallery modes provide examples of spectral
clusters which concentrate in a $\lambda^{-\frac{2}{3}} \times
\lambda^{-\frac{n-2}{2}}$ neighborhood of a boundary geodesic (see
Grieser~\cite{grieser}).  These examples show that one cannot
achieve linear spectral cluster estimates better than
\begin{equation}\label{E:linest}
\|\chi_\lambda\|_{L^2 \to L^q} \lesssim \begin{cases}
\lambda^{(\frac{2}{3}+\frac{n-2}{2})(\hf-\frac{1}{q})} &
2 \leq q \leq \frac{6n+4}{3n-4}\\
\lambda^{n(\hf-\frac{1}{q})-\hf} &  \frac{6n+4}{3n-4} \leq q \leq
\infty
\end{cases}
\end{equation}
The estimates \eqref{E:linest}
were recently proven for dimension $n=2$ in~\cite{SmSoBdry},
along with partial results in higher dimensions.  The
question of whether or not they hold in general
in higher dimensions remains open.

In this work, we establish the following multilinear spectral
cluster estimates on a general Riemannian manifold with boundary.
We restrict attention to dimensions $n=2,3$, since that
is where our results are in some context sharp.

\begin{theorem}\label{T:mainthm}
Let $(M^n,\g)$ and $\chi_\lambda$ be as above, with either
Dirichlet or Neumann eigenfunctions, and
let $\l\ge\mu\ge\nu$.  Then the
following bilinear estimate holds 
\begin{equation}\label{E:bilin}
\| \chi_\l f \,\chi_\mu g\|_{L^2(M)} \lesssim \Lambda(\mu)\,
 \|f\|_{L^2(M)} \|g\|_{L^2(M)}
\end{equation}
with $\Lambda$ defined as
\begin{equation*}
\Lambda(\mu) = \begin{cases} \mu^\frac{1}{3} & n=2\\
\mu^{\frac{2}{3}}(\log{\mu})^\hf & n=3
\end{cases}
\end{equation*}
In addition, the following trilinear estimate
holds for $n=2,3$
\begin{equation}\label{E:multilin}
\| \chi_\l f \,\chi_\mu g \,\chi_\nu h\|_{L^2(M)} \lesssim
(\mu\nu)^{\frac{3n-4}{6}}\|f\|_{L^2(M)}
\|g\|_{L^2(M)} \|h\|_{L^2(M)}
\end{equation}
\end{theorem}

For $n=2$, the bilinear and trilinear estimates imply
the estimate \eqref{E:linest}, respectively for $q=4$ and $q=6$.
Moreover, \eqref{E:multilin} implies \eqref{E:bilin}
for Neumann conditions by taking $h$ constant.
For $n=3$ this is no longer the case. However, the bilinear estimate
for $n=3$ implies (up to a logarithmic loss) 
the best possible $L^4$ linear estimate for
manifolds with Lipschitz metric, by the examples of \cite{SmSolowreg},
and our proof in fact establishes Theorem \ref{T:mainthm} in this
context.

As in \cite{SmSoBdry}, a key step is to work on the double $\tilde M$
of the manifold $M$, obtained by attaching two copies of $M$ along
the boundary, and taking coordinate patches along 
$\partial M\subset \tilde M$ which agree with
geodesic normal coordinates $(x',x_n)$
on each copy of $M$. 
In these coordinates, the lift $\tilde \g$ of $\g$ to $\tilde M$ is
given by $\g^{ij}(x',|x_n|)$, and hence $\tilde \g$ extends to
$\tilde M$ with a Lipschitz type
singularity along $\partial M$. Dirichlet and Neumann eigenfunctions
on $M$ correspond to eigenfunctions on $\tilde M$ which
are, respectively, odd or even under $x_n\to -x_n$. Hence, $L^q$
bounds on spectral clusters, and the multilinear analogs we consider,
can be obtained by proving the same bounds on $(\tilde M,\tilde\g)$.

The linear estimates of \cite{SmSoBdry} were obtained by establishing
mixed-norm $L^q_xL^2_t$ estimates on $\tilde M$ for
the evolution of a spectral cluster
under the wave equation, and we follow a similar approach here.
In that paper, the precise nature of the singularity of $\tilde \g$ along
$\partial M$ was used, and a microlocal decomposition of the cluster
was made in terms of angle from tangent to $\partial M$.
Estimates were obtained over small slabs, with size
depending on the frequency and angle, and summing over slabs led
to a frequency dependent loss for the estimates on $\tilde M$.

In contrast, the results of this paper go through generally for the case of a
boundary-free Riemmanian manifold with metric of Lipschitz regularity,
as with the linear spectral cluster estimates of \cite{SmC1alph}, or 
the Strichartz estimates of Tataru \cite{Tat}. As in those papers,
we obtain estimates over small slabs with size depending on the
frequency, and use a rescaling argument to reduce matters to
obtaining estimates for $C^2$ metrics. 
We then use wave packet methods to obtain dispersive estimates,
as in \cite{KT}, \cite{SmC2} and \cite{Tat}, 
Summing over slabs then leads to a frequency dependent loss.

\section{Microlocal Reductions}
For the remainder of this paper, we assume $M$ is a compact manifold
without boundary, and $\g$ is a metric of Lipschitz regularity.
The condition that $f,g,h$ be spectrally localized
can be relaxed, and we work instead with a {\em quasimode} condition.
We state the condition
for $f$ here, the condition for $g$ and $h$ being analogous.
For each local coordinate chart we write
\begin{equation}\label{quasimode}
\g\,\dsq f+\l^2 f=w\,,\qquad \g\,\dsq=
\sum_{j,k} \g^{jk}(x)\,\partial_j\partial_k\,.
\end{equation}
Given such an equation, we set
$$
|||f|||_\l=
\|f\|_{L^2}+\l^{-1}\|df\|_{L^2}+\l^{-2}\|\dsq f\|_{L^2}+
\l^{-1}\|w\|_{L^2}\,.
$$

If $f=\chi_{\l}f$, then for $\phi$ a smooth cutoff to
a local coordinate system, the function
$\phi f$ satisfies the equation \eqref{quasimode} on $\R^n$,
and 
$$
|||\phi f|||_\l\lesssim \|f\|_{L^2(M)}\,.
$$
For the bilinear estimates it thus suffices to prove,
for each coordinate chart, that
$$
\|\phi f\,\phi g\|_{L^2}\lesssim 
\Lambda(\mu)\,|||\phi f|||_\l\,|||\phi g|||_\mu\,,
$$
and analogously for the trilinear version.

By choosing appropriate coordinates, we
may assume $f$ satisfies \eqref{quasimode} on $\R^n$, with $f$ supported
in the unit ball, and
$$
\|\g^{ij}-\delta^{ij}\|_{\mathrm{Lip}(\R^n)}\le c_0\,,
$$
with $c_0$ a constant to be chosen suitably small.

Let $S_r=S_r(D)$ denote a smooth cutoff on the Fourier transform side to
frequencies of size $|\xi|\le r$.
Let $\g_\l=S_{c^2\l}\g$, for $c$ to be chosen suitably small.
Then
$$
\|(\g-\g_\l)\dsq f\|_{L^2}\lesssim c^{-2}\l^{-1}\|\dsq f\|_{L^2}\,,
$$ 
and thus we may replace $\g$ by $\g_\l$ in \eqref{quasimode}
at the expense of absorbing the above
term into $w$, which does not change the size of $|||f|||_\l$.

We next take a microlocal partition of unity, $1=\sum \Gamma(D)$, where
each $\Gamma(\xi)$ is a smooth symbol of order $0$ supported in a cone of
small angle. By the Coifman-Meyer commutator theorem \cite{CM}, since
$\g$ is Lipschitz
\[
[\g,\Gamma(D)]\,d\,:\,L^2(\R^n)\rightarrow L^2(\R^n)\,,
\]
hence $\Gamma(D)f$ satisfies the equation \eqref{quasimode}, with
$|||\Gamma(D)f|||_\l\lesssim |||f|||_\l\,.$

Since there are finitely many terms, 
we may replace $f$ by $\Gamma(D)f$, which is 
no longer compactly supported,
but is rapidly decreasing and smooth outside the unit ball.
Without loss of generality
we assume that $\Gamma(\xi)$ is supported within a small angle 
of the $\xi_1$ axis. 
We similarly replace $g$ and $h$ by $\Gamma'(D)g$ and $\Gamma''(D)h$,
localized in frequency to small cones along general directions.

Letting $x'=x_2$ in case of dimension $n=2$, and 
$x'=(x_2,x_3)$ in case of dimension $n=3$, we bound
\begin{align*}
\|fg\|_{L^2}&\le 
\|f\|_{L^\infty_{x_1}\!L^2_{x'}}\,
\|g\|_{L^2_{x_1}\!L^\infty_{x'}}\\
\\
\|fgh\|_{L^2}&\le 
\|f\|_{L^\infty_{x_1}\!L^2_{x'}}\,
\|g\|_{L^4_{x_1}\!L^\infty_{x'}}\,
\|h\|_{L^4_{x_1}\!L^\infty_{x'}}
\end{align*}
Since $\Gamma(D)f$ is rapidly decreasing outside the unit ball,
it suffices to take the norms above over the ball of radius 2.
Theorem \ref{T:mainthm} is then a result of the following.

\begin{theorem}
Suppose that $f$ satisfies the equation
\begin{equation}\label{lquasimode}
\g_\l \dsq f+\l^2 f=w
\end{equation}
Then the following hold, where the norms on the left side are
over a bounded set
\begin{align}\label{est1}
\|f\|_{L^2_{x_1}\!L^\infty_{x'}}&\lesssim 
\l^{\frac 23}(\log \l)^{\frac 12}\,|||f|||_\l\,,\qquad n=3
\\\nonumber
\\\label{est2}
\|f\|_{L^4_{x_1}\!L^\infty_{x'}}&\lesssim
\begin{cases} 
\l^{\frac 13}\,|||f|||_\l\,,\qquad n=2\\
\l^{\frac 56}\,|||f|||_\l\,,\qquad n=3
\end{cases}
\end{align}
Futhermore, if $\hat f(\xi)$ is supported in a small cone about the $\xi_1$
axis, then
\begin{equation}\label{est3}
\|f\|_{L^\infty_{x_1}\!L^2_{x'}}\lesssim |||f|||_\l
\end{equation}
\end{theorem}
\noindent{\it Proof.}
We start by localizing $f$ dyadically in frequency. Let
$$
f=S_{c\l}f+ (S_{c^{-1}\l}-S_{c\l})f +(1-S_{c^{-1}\l})f
\equiv
f_{<\l}+f_\l+f_{>\l}\,.
$$
Since $[S_{c\l},\g_\l]\,d\,:\,L^2\rightarrow L^2$, $f_\l$
satisfies \eqref{lquasimode}, with $|||w_\l|||_\l\lesssim |||f|||_\l\,,$
and similarly for $f_{>\l}$ and $f_{<\l}$. Furthermore, by the frequency
localization of $\g_\l$, each of $w_\l$, $w_{<\l}$, and $w_{>\l}$ is 
also localized to the appropriate range of frequencies.

A simple integration by parts argument (see the proof of Corollary 5 of 
\cite{SmC1alph})
yields that, for $c$ sufficiently small,
\[
\l\|f_{<\l}\|_{L^2(\R^n)}+\|d f_{>\l}\|_{L^2(\R^n)}\lesssim |||f|||_\l\,.
\]
Elliptic regularity additionally gives the bound
\[
\|\dsq f_{>\l}\|_{L^2(\R^n)}\lesssim \l\,|||f|||_\l\,.
\]
Sobolev embedding then yields each of the estimates \eqref{est1}--\eqref{est3}
for $f_{<\l}$ and $f_{>\l}$. Indeed, there is a gain of 
$\l^{\frac 23}(\log\l)^{\frac 12}$ 
in the estimate \eqref{est1}, and a gain of $\l^{\frac 7{12}}$ in 
the estimate \eqref{est2},
for these terms. 

Consequently, we are reduced to establishing
\eqref{est1}--\eqref{est3} for the term $f_\l$.
We start with \eqref{est3}.
Let $V$ denote the vector field
$$
V=2(\partial_1\!f_\l)\,\g_\l\,df_\l+
\bigl(\l^2 f_\l^2-\langle \g_\l\,df_\l,df_\l\rangle\bigr)
\overrightarrow{e_1}\,.
$$
Then
$$
\mathrm{div}\, V=2(\partial_1\!f_\l)\,(\mathrm{div}\,\g_\l)\cdot df_\l
+2(\partial_1\!f_\l)\,w_\l-\langle (\partial_1 \g_\l)df_\l,df_\l\rangle\,.
$$
Applying the divergence theorem on the set $x_1\le r$ yields
$$
\int_{x_1=r}V_1\,dx'
\lesssim
\l^2\|f_\l\|_{L^2(\R^n)}^2+\|df_\l\|_{L^2(\R^n)}^2+\|w_\l\|_{L^2(\R^n)}^2\,.
$$
Since $\g_\l$ is pointwise close to the flat metric,
we have pointwise that
$$
V_1\ge \tfrac 34 |\partial_1\! f_\l|^2+\tfrac 34 \l^2|f_\l|^2
-|\partial_{x'}\!f_\l|^2\,.
$$
The frequency localization of $\widehat f_\l$ to $|\xi'|\le c\l$ yields
$$
\int_{x_1=r} V_1\,dx'\ge
\frac 12\int_{x_1=r}|df_\l|^2+\l^2|f_\l|^2\,dx'\,.
$$
Consequently,
\begin{align*}
\l^{-1}\|df_\l\|_{L^\infty_{x_1}L^2_{x'}}+
\|f_\l\|_{L^\infty_{x_1}L^2_{x'}}&\lesssim
\|f_\l\|_{L^2(\R^n)}+\l^{-1}\|df_\l\|_{L^2(\R^n)}+\l^{-1}\|w_\l\|_{L^2(\R^n)}
\\
&\le |||f_\l|||_\l\,,
\end{align*}
yielding estimate \eqref{est3}.

For the remainder, we assume that
$\hat f(\xi)$ is localized to a small cone along the direction $\omega$.
In this case, the above argument yields uniform $L^2$
bounds over hyperplanes of the
form $\omega\cdot x=r$. In the proof of \eqref{est1}--\eqref{est2}
we will use the following consequence. Suppose
that $S_R$ is a slab of the form $\omega\cdot x\in I$, where $I$ is an
interval of length $|I|=R$. Then
\begin{equation}\label{fluxest}
\l^{-1}\|df_\l\|_{L^2(S_R)}+
\|f_\l\|_{L^2(S_R)}\lesssim
R^{\frac 12}|||f_\l|||_\l\,.
\end{equation}

We cover the bounded set on which the norms in \eqref{est1} and \eqref{est2}
are taken by $\approx \l^{\frac 13}$ slabs of the form $\omega\cdot x\in I$,
where $|I|=R=\l^{-\frac 13}$. Then
\begin{equation}\label{sumbound}
\|f_\l\|_{L^p_{x_1}\!L^\infty_{x'}}\lesssim  
\l^{\frac 1{3p}}\sup_{S_R} \|f_\l\|_{L^p_{x_1}\!L^\infty_{x'}(S_R)}\,.
\end{equation}
We will establish the following result. Suppose that $Q_R$ is a cube
of sidelength $R=\l^{-\frac 13}$, and $Q_R^*$ its double. Then
\begin{multline}\label{est'}
\|f_\l\|_{L^p_{x_1}\!L^\infty_{x'}(Q_R)}\\
\lesssim 
c_p(\l)R^{-\frac 12}
\bigl(\,\|f_\l\|_{L^2(Q_R^*)}+\l^{-1}\|df_\l\|_{L^2(Q_R^*)}
+R\l^{-1}\|w_\l\|_{L^2(Q_R^*)}\bigr)
\end{multline}
where
\begin{equation*}
c_2(\l)=\l^{\frac 12}(\log\l)^{\frac 12}\,,\quad n=3
\qquad\qquad
c_4(\l)=
\begin{cases}
\l^{\frac 14}\,,\quad n=2\\
\l^{\frac 34}\,,\quad n=3
\end{cases}
\end{equation*}
If we cover the slab $S_R$ by disjoint cubes $Q_R$, then
by \eqref{fluxest} we obtain
\begin{equation*}
\|f_\l\|_{L^p_{x_1}\!L^\infty_{x'}(S_R)}\lesssim 
c_p(\l)\,|||f_\l|||_\l\,,
\end{equation*}
and \eqref{sumbound} yields \eqref{est1}--\eqref{est2}.

The estimate \eqref{est'} is scale-invariant. Precisely, if we change
$x\rightarrow Rx$, so that $Q$ becomes a cube of size 1,
and $f_{\l}(R\,\cdot)$ is frequency localized at scale $R\l=\l^{\frac 23}$,
then, with $\mu=\l^{\frac 23}$, estimate \eqref{est'} is equivalent to
the following 
\begin{equation}\label{est''}
\|f_\mu\|_{L^p_{x_1}\!L^\infty_{x'}(Q)}\\
\lesssim 
c_p(\mu)
\bigl(\,\|f_\mu\|_{L^2(Q^*)}+\mu^{-1}\|df_\mu\|_{L^2(Q^*)}
+\mu^{-1}\|w_\mu\|_{L^2(Q^*)}\bigr)\,.
\end{equation}
Here, $f_\mu(x)=f_\l(\l^{-\frac 13}x)$, which satisfies the equation
$$
\g_\mu\dsq f_\mu+\mu^2 f_\mu=w_\mu\,,
$$
with $\g_\mu(x)=\g_\l(\l^{-\frac 13}x)$.
Observe that
\begin{equation}\label{C2est1}
\|d\g_\mu\|_{L^\infty}\le c_0\l^{-\frac 13}=c_0\mu^{-\frac 12}\,,
\end{equation}
hence
$$
\|\g_\mu-\g_{\mu^{1/2}}\|_{L^\infty}\le c_0\mu^{-1}\,,
$$
with $\g_{\mu^{1/2}}=S_{c^2\mu^{1/2}}\g_\mu\,.$ Thus
$f_\mu$ satisfies the equation
\begin{equation}\label{C2eqn}
\g_{\mu^{1/2}}\dsq f_\mu+\mu^2 f_\mu=w_\mu\,,
\end{equation}
with the right-hand side of \eqref{est''} of comparable size.

It follows from \eqref{C2est1} that the metric $\g_{\mu^{1/2}}$ is
of regularity $C^2$. Indeed,
$$
\|\g^{jk}_{\mu^{1/2}}-\delta^{jk}\|_{C^2}\le c_0\,.
$$
To establish \eqref{est''},
we may thus use techniques developed to establish dispersive
estimates for 
operators of principal type with $C^2$ coefficients. 
We follow below the path through squarefunction estimates for solutions
to a first order hyperbolic equation, by introducing a time
variable, as in \cite{SmC2}.
It should also be possible to establish the dispersive
estimates directly for
\eqref{C2eqn} using methods of \cite{KT}.

Let
$$
p(\cd,\xi)=
S_{c^2\mu^{1/2}}
\Bigl(\,\sum_{j,k}\g_{\mu^{1/2}}^{jk}(\cdot)\,\xi_j\,\xi_k\Bigr)^{\frac 12}\,.
$$
Then
$$
\|p(x,D)^2 f_\mu+\g_{\mu^{1/2}} \dsq f_\mu\|_{L^2(\R^n)}\lesssim
\mu\,\|f_\mu\|_{L^2}\,.
$$
Thus,
$$
\bigl(\mu+p(x,D)\bigr)\bigl(\mu-p(x,D)\bigr)f_\mu=w_\mu\,,
$$
with the error harmlessly absorbed into $w_\mu$.
The operator $\mu+p(x,D)$ is elliptic on the frequency support
of $\bigl(\mu-p(x,D)\bigr)f_\mu$, hence we may write
$$
\bigl(\mu-p(x,D)\bigr)f_\mu=\mu^{-1}w_\mu\,,
$$
with $|||f|||_\l$ still of comparable size.
Finally, let
$$
u(t,x)=e^{-it\mu}f_\mu(x)\,,\qquad\qquad F=\mu^{-1}e^{-it\mu}w_\mu\,.
$$
Then
$$
\bigl(\partial_t+ip(x,D)\bigr)u=F\,,
$$
and it suffices to show that
\begin{equation}\label{dispest}
\|u\|_{L^p_{x_1}\!L^\infty_{x'}\!L^2_t(Q\times[0,1])}
\lesssim
c_p(\mu)\bigl(\,
\|u\|_{L^\infty_tL^2_x([0,1]\times\R^n)}+
\|F\|_{L^1_tL^2_x([0,1]\times\R^n)}\bigr)\,.
\end{equation}

Our proof of \eqref{dispest} follows very closely the proof of
the linear spectral cluster estimates in \cite{SmC2}; we
outline just the main steps here.
Following \cite[\S 3]{SmC2}, consider
the wave packet transform of $u_\mu$,
$$
\bigl(T_\mu u\bigr)(t,x,\xi)=\mu^{\frac n4}\,
\int e^{-i\langle \xi,z-x\rangle}\,\phi\bigl(\mu^{\frac 12}(z-x)\bigr)
\,u(t,z)\,dz\,,
$$
where $\phi$ is a real, even Schwartz function, with 
$\|\phi\|_{L^2}=(2\pi)^{-\frac n2}$, 
and with Fourier transform
supported in the unit ball $\{|\xi|\le 1\}\,.$
Then
$$
\partial_t T_\mu u(t,x,\xi)=
\Bigl(d_\xi p(x,\xi)\cdot d_x-d_x p(x,\xi)\cdot d_\xi\Bigr)
T_\mu u(t,x,\xi)
+G(t,x,\xi)\,,
$$
where $G(t,x,\xi)=0$ unless $\frac 18\mu<|\xi|<2\mu$, and
$$
\|G\|_{L^1_tL^2_{x,\xi}}\lesssim 
\|u\|_{L^1_tL^2_x}+\|F\|_{L^1_tL^2_x}\,.
$$
Let $\chi_{t}$ denote the canonical transform on 
$\R^{2n}_{x,\xi}=T^*(\R^n)$ generated by the Hamiltonian flow of $p$.
Thus, $\chi_{t}(x,\xi)=\gamma(t)$, where $\gamma$ is the integral
curve with $\gamma(0)=(x,\xi)$. Then we have
$$
\bigl(T_\mu u\bigr)(t,x,\xi)=\bigl(T_\mu u\bigr)(0,\chi_{-t}(x,\xi))+
\int_0^t G(r,\chi_{r-t}(x,\xi))\,dr\,.
$$
Thus, $T_\mu u(t,x,\xi)$ is an integrable superposition
over $r$ of $1_{t>r}$ multiplied by a 
function invariant under the Hamiltonian flow of $p$. 

Since $u(t,x)=T_\mu^* \bigl(T_\mu u\bigr)(t,x,\xi)$, it suffices to show
\begin{equation}
\|W\tilde f\|_{L^p_{x_1}\!L^\infty_{x'}\!L^2_t(\R^n\times[0,1])}
\lesssim
c_p(\mu)\|f\|_{L^2_{x,\xi}}\,,
\end{equation}
where
$$
\bigl(W\!\tilde f\,\bigr)(t,x)=T_\mu^*\bigl(\tilde f\circ\chi_{-t}\bigr)(x)\,.
$$
This is in turn equivalent to the following bounds
\begin{equation}\label{dispest'}
\|WW^*\!F\|_{L^p_{x_1}\!L^\infty_{x'}\!L^2_t(Q\times[0,1])}\lesssim
c_p(\mu)^2\,
\|F\|_{L^{p'}_{x_1}\!L^1_{x'}\!L^2_t(Q\times[0,1])}\,.
\end{equation}

The operator $WW^*$ has an integral kernel $K$ which
is highly localized
to a $\mu^{-1}$ neighborhood of the light cone, with the dispersive
rate of decay away from the origin, see \cite[(3.11)]{SmC2},
\begin{equation}\label{Kest}
|K(s,y;t,z)|\lesssim
\mu^n\bigl(\,1+\mu\,|y_1-z_1|\,\bigr)^{-\frac{n-1}2}
\,\bigl(\,1+\mu\,\bigl|\,d(y,z)-|s-t|\bigr|\,\bigr)^{-N}\,,
\end{equation}
with $d(y,z)$ the distance of $y$ to $z$ determined by $p$.

We remark that in \cite{SmC2} this estimate
was established assuming the kernel was microlocalized near
the $\xi_1$ axis. That assumption, however, was necessary
for $L^2$-energy estimates, not the above dispersive estimates.
Indeed, the proof of \cite[(3.11)]{SmC2}
establishes \eqref{Kest} with $|y_1-z_1|$ replaced by $|y-z|$,
since $|t-s|\approx |y-z|$ on the light cone,
and hence holds without any assumption of conic micro-localization.

Estimate \eqref{Kest} implies that, for each $(y_1,z_1)$, 
$$
\Bigl\|\int K(s,y;t,z)\,v(t,z')\,dt\,dz'\Bigr\|_{L^\infty_{y'}L^2_s}
\lesssim
\mu^{n-1}(1+\mu\,|y_1-z_1|)^{-\frac{n-1}2}
\|v\|_{L^1_{z'}L^2_t}\,.
$$
For $n=3$ and $p=2$, estimate \eqref{dispest'} follows from
$$
\int_{|z_1|\le 2}
\mu^2(\,1+\mu\,|y_1-z_1|\,)^{-1}\le \mu\log\mu\,.
$$
For $p=4$, \eqref{dispest'}
follows from the Hardy-Litlewood-Sobolev inequality,
together with the bound
$$
\qquad\qquad\mu^{n-1}(\,1+\mu\,|y_1-z_1|\,)^{-\frac{n-1}2}
\le \mu^{n-\frac 32}|y_1-z_1|^{-\frac 12}\,.\qquad\qquad\qed
$$

\end{document}